\documentclass[12pt]{amsart}
   \newtheorem{lemma}{Lemma}[section]
   \newtheorem{theorem}[lemma]{Theorem}
   \newtheorem{remark}[lemma]{Remark}
   \newtheorem{prop}[lemma]{Proposition}
   \newtheorem{coro}[lemma]{Corollary}

\newcommand{\Om}{\Omega}

\newcommand{\p}{\partial}
\renewcommand{\phi}{\varphi}

\newcommand{\R}{{\mathbb R}}

\title[Stochastic Swift-Hohenberg System]
{Ergodic Dynamics of the Stochastic Swift-Hohenberg System}

\thanks{This work was partly supported by the NSF Grants DMS-0209326
and DMS-0112351, the National Natural Science Foundation of China
Grant 10171044, the Natural Science Foundation of Jiangsu Province
Grant BK2001024, and the Foundation for University Key Teachers of
the Ministry of Education of China. }

\author{Wei Wang }
\address[W.~Wang ]
{Department of Mathematics\\
Nanjing University \\
Nanjing 210008, China }
\email[W.~Wang]{wangshengtaoo@sina.com}

\author{Jianhua  Sun }
\address[J.~Sun ]
{Department of Mathematics\\
Nanjing University \\
Nanjing 210008, China }
\email[J.~Sun]{jhsun@nju.edu.cn}

\author{Jinqiao Duan}
\address[J.~Duan]
{Department of Applied Mathematics\\
Illinois Institute of Technology\\
Chicago, IL 60616, USA}
\email[J.~Duan]{duan@iit.edu}

\date{August 18, 2004(final version)}


\begin{document}


\begin{abstract}

{\bf Nonlinear Analysis B}: in press, 2004.

 The Swift-Hohenberg fluid convection system  with both local
and nonlocal nonlinearities  under the influence of   white noise
is studied. The objective is to understand the difference in the
dynamical behavior in both local and nonlocal cases. It is proved
that when sufficiently many of its Fourier modes are forced, the
system has a unique invariant measure, or equivalently, the
dynamics is ergodic.  Moreover, it is found that the number of
modes to be stochastically excited for ensuring the ergodicity in
the local Swift-Hohenberg system depends {\em only} on the
Rayleigh number (i.e., it does not even depend on the random term
itself), while this number for the nonlocal Swift-Hohenberg system
relies additionally on the bound of the kernel in the nonlocal
interaction (integral) term, and on the random term itself.

\bigskip
Submitted to:  {\bf Nonlinear Analysis B}
\end{abstract}

\maketitle

\section{Introduction}

Density gradient-driven fluid convection  arises in geophysical
fluid flows in the atmosphere, oceans and the earth's mantle. The
Rayleigh-Benard convection is a prototypical model for fluid
convection, aiming at predicting spatio-temporal convection
patterns. The mathematical model for the Rayleigh-Benard
convection involves the Navier-Stokes equations coupled with the
transport equation for temperature. When the Rayleigh number is
near the onset of the convection, the Rayleigh-Benard convection
model may be approximately reduced to an amplitude or
order parameter equation, as derived by Swift-Hohenberg \cite{SH77}.\\

In the literature, most works (e.g.,\cite{HMBD95,GL95,MS95}) on
the Swift-Hohenberg model deal with the following  evolution
equation for order parameter $u(x,t)$, which is localized  version
of the model originally derived by Swift-Hohenberg \cite{SH77},
\begin{equation*}
u_t=\varrho u-(1+\partial_{xx})^2u-u^3,
\end{equation*}
where $\varrho$ measures the difference of the Rayleigh number
from its critical onset value and the cubic term $u^3$ is an yet
``approximation" of a nonlocal integral term in the
original Swift-Hohenberg model \cite{SH77}.\\

Roberts (\cite{Rob92,Rob95}) re-examined the rationale for using
the Swift-Hohenberg model as a reliable simplified model of the
spatial pattern evolution in fluid convection. He argued that,
although the localization approximation used in (1) makes some
sense, the approximation is deficient in describing some basic
features of such systems, and devised via symmetry argument the
following modified Swift-Hohenberg equation with nonlocal
interactions
\begin{equation*}
u_{t}=\varrho u-(1+\p_{xx})^2u-uG*u^2,
\end{equation*}
where $G*u^2$ is a spatial convolution integral and $G(\cdot)$ is
a given radially symmetric positive (or nonnegative) function. We
call this function $G$ the kernel for the nonlocal term. In fact,
nonlocal integral terms often appear naturally in amplitude
equation models for nonequilibrium systems; see, e.g.,
\cite{MV, nonlocal_ks, Duan-Ervin,LGDE00}.\\

Our goal in this paper is to examine the above local and nonlocal
Swift-Hohenberg models, by investigating the dynamical difference
of both models under random impact as well as under nonlocal
interactions.\\

Fluid systems are often subject to random environmental
influences. On the one hand, there is a growing recognition of a
role for the inclusion of stochastic effects in the modeling of
complex systems. Randomness can have delicate impact on the
overall evolution of such systems, for example, noise-induced
phase transition or stochastic bifurcation \cite{CarLanRob01},
stochastic resonance \cite{Imkeller},
and noise-induced pattern formation \cite{Gar,BlomkerStani}.
Taking stochastic effects into account is of central importance
for the development of mathematical models of such complex
phenomena in engineering and science. Macroscopic models in the
form of partial differential equations for these systems contain
such randomness as stochastic forcing, uncertain parameters,
random sources or inputs, and random   boundary conditions.
Stochastic partial differential equations (SPDEs) are
appropriate models for randomly influenced spatially extended
systems. On the other hand, the inclusion of such effects has
led to interesting new mathematical problems  at the interface
of probability  and partial differential equations.\\

In this paper, we consider the Swift-Hohenberg model, taking
stochastic forcing as well as nonlocal interactions, into
account.\\

First, we consider the local stochastic Swift-Honenberg system
({\bf LSSH})
\begin{equation}  \label{LSSH}
u_t=\varrho u-(1+\partial_{xx})^2u-u^3+F(x,t),
\;\;u(0)=u_0,
\end{equation}
under the periodic boundary condition on with period $2\pi$, in
which the stochastic force is taken to be of the form
\begin{equation}\label{Noise1}
F(x,t)=\partial_tW(x,t)
      =\partial_t\sum_{i=1}^{N'}b_ie_i(x)w_i(t),
      \;\; N'\leq \infty,
\end{equation}
where $\{e_{i}\}$ is an orthonormal basis of the Hilbert space
$H=L^{2}_{per}(0, 2\pi)$ formed by the eigenvectors of the
operator $A=-(1+\partial_{xx})^{2}$, $\{w_{i}(t), t\geq 0\}$ are
independent standard Wiener processes, and the real coefficients
$b_{i}\geq 0$ satisfies
\begin{equation*}
b_{i}\neq 0\;,\;\;1\leq i\leq N\leq N'\label{Mode}
\end{equation*}
for some sufficiently large $N$.\\

Then, we study the nonlocal stochastic Swift-Hohenberg system
({\bf NLSSH})
\begin{equation}\label{NLSSH}
u_t=\varrho u-(1+\partial_{xx})^2u-uG*u^2+\tilde{F}(x,t),
\;\;u(0)=u_0,
\end{equation}
with a positive kernel (i.e., $G>0$) and a special nonnegative
kernel (i.e., $G\geq 0$), respectively, under the periodic
boundary condition on $I$, in which
\begin{equation*}
 G*u^2=\int_0^{2\pi}G(x-\xi)u^2(\xi,t)d\xi,
\end{equation*}
and the stochastic force is taken to be of the form
\begin{equation}\label{Noise2}
\tilde{F}(x,t)=\partial_t\tilde{W}(x,t)
      =\partial_t\sum_{i=1}^{N'}\tilde{b}_ie_i(x)w_i(t),
      \;\; N'\leq \infty.
\end{equation}

The diagnostic tool we choose to compare the dynamical
differences between the local and nonlocal stochastic
Swift-Hohenberg systems is {\em ergodicity}. Namely, we are
interested in the ergodicity of the Markov solution process in $H$
generated by the stochastic Swift-Hohenberg equation (\ref{LSSH})
which forms a random dynamical system (RDS). More precisely, we
first prove that if sufficiently many of its Fourier modes are
forced, the LSSH system has a unique invariant measure, or
equivalently, the dynamics is ergodic in the phase space, and that
all solutions converge exponentially fast in distribution to this
unique measure ({\bf Theorem \ref{contraction} and Corollary
\ref{ergodicity}}). Our results show that under certain conditions
the RDS generated by (\ref{LSSH}) has a global exponentially
attracting fixed point in the sense of distribution. And this
fixed point is the unique invariant measure which is supported on
the random attractor. Namely, the LSSH system (\ref{LSSH}) is
ergodic. To prove our results, we decompose the LSSH system
(\ref{LSSH}) into a low mode part with a finite dimension and a
high mode part with an infinite dimension. In the low mode part,
we use the maximal coupling approach, which is used in
\cite{KS03}, to prove that the low mode part can be coupled. The
high mode part can be enslaved, and then
the ergodic result is obtained; see also \cite{Hair02,Mat02}.\\

We   verify that the dynamics of the nonlocal stochastic
Swift-Hohenberg equation (\ref{NLSSH}) is also ergodic, provided
sufficiently many of its Fourier modes are forced ({\bf Theorems
\ref{im} and \ref{im2}}). However, our interest here is to
investigate the difference in the conditions ensuring the
ergodicity of the LSSH and NLSSH systems. We find out that (see
discussions in Sec. 7):

(i) For the {\em local} Swift-Hohenberg system, the estimated
number of Fourier modes to be randomly excited for ensuring
ergodicity depends only on the parameter $\varrho$, which
measures the difference of the Rayleigh number from its critical
onset value. Note that this number does not depend on the random
forcing term;

(ii) For the {\em nonlocal} Swift-Hohenberg system  with positive
kernel $G$ in the nonlocal nonlinearity, the  estimated number of
Fourier modes to be randomly excited for ensuring ergodicity
depends additionally on the upper and lower bounds of the positive
kernel, and on the random term itself.  \\

Recently, there have been a number of papers on ergodicity of
stochastically-forced partial differential equations (SPDEs) and,
in particular, 2D Navier-Stokes equations. The stochastic force of
the SPDEs may be in the form of \ref{Noise1} above or in  the
following  form (so-called  kick-force)
\begin{equation}
F(x,t)=\sum_{k=-\infty}^{\infty}f_{k}\delta(t-Tk)
\;,\;\;f_{k}=\sum_{j=1}^{\infty}b_{j}\xi_{jk}e_{j}
\end{equation}
where $b_{j}\geq 0$ are some constants such that $\sum
b_{j}^{2}<\infty$ and $\{\xi_{jk}\}$ are independent random
variables with $k$-independent distributions, or the white (in
time) force of the form (\ref{Noise1}). A coupling approach and
the so-called Kantorovich functionals (e.g., \cite{KS01,
Kuk02,Mat02}) have been developed  to show exponential convergence
to the unique invariant measure for stochastic Novier-Stokes equations.
Similar results were obtained in \cite{MY02}. For the white noise
forced PDE, \cite{FM95} and \cite{Mat99} obtained the ergodicity
of the stochastic forcing Navier-Stokes equation for the case when
the random force is singular in $x$. For the two-dimensional
Navier-Stokes equation with periodic boundary condition and random
forcing, \cite{EMS01} proved uniqueness of the stationary measure
under the condition that all ``deterministic modes'' are forced,
by studying the Gibbsian dynamics of the low modes. The idea in
\cite{EMS01} has been used in other papers to get explicit
results. Moreover, \cite{Hair02} obtained the exponential mixing
property of the stochastic Ginzburg-Landau equation by the
so-called $binding$ $construction$. We also mention that
\cite{KS03} developed the idea in \cite{KS01,Kuk02} to get the
rates of the measure converge to the unique invariant measure for
2D Navier-Stokes equation with a white noise to all Fourier modes.\\

In the present paper, the linear part $-(1+\partial_{xx})^2 u$ in
the Swift-Hohenberg models is not dissipative, i.e., it has
positive eigenvalues. Similar situation has also been considered
in \cite{Hair02}, for example.  In this case, the coupling
approach can still be used to investigate ergodicity. In fact, we
can directly study the maximal coupling solutions of the
stochastic Swift-Hohenberg system, without constructing the
Kantorovich functional; see \cite{Hair02,Mat02}.   Then the key
fact that the orbit in set $S(m,k)$ comes from $R(m)\cup S(m,m)$
helps us to prove the main result without using the exponential
decay of a Kantorovich type functional, for details see Sec.4 and
Sec.5. For simplicity of presentation, in this paper we only work
with the periodic boundary condition and only consider the case
when only finite modes are forced, i.e, $N'<\infty$. However, we
emphasize that our argument applies in principle to the case when
the other more physical boundary conditions are posed and when all
the modes are subject to random excitation. Moreover, our
techniques also applies to more realistic two-dimensional
Swift-Hohenberg model, i.e., $(1+\p_{xx})^2$ is replaced by
$(1+\p_{xx} + \p_{yy})^2$. \\

We remark that a systematic treatment for random dynamical
systems (including ergodicity) is emerging; see, e.g.,
\cite{Arn98,PZ96,
CF94, WS02, WSD02}. \\

The remainder of this paper is organized as follow. In Sec. 2, we
decompose the local system LSSH   into the low mode part and the
high mode part. Sec. 3 is devoted to some energy estimates of the
solution. The coupling construction and a coupling result are
given in Sec. 4. Then we give the ergodic result in Sec. 5. In
Sec. 6, we study the ergodicity of nonlocal system NLSSH with
positive or nonnegetive kernels of nonlocal nonlinearity,
respectively, by the techniques developed in the previous
sections. Finally, we discuss the difference between LSSH system
and NLSSH system by analysing the minimal conditions ensuring the
ergodicity in Sec. 7.


\section{Decomposition of the Local System }

We consider the following system on the real line
\begin{eqnarray}
u_t-\varrho u+(1+\partial_{xx})^2u+u^3&=&F(x,t)\;,\label{LSSH1}\\
u(0) &=& u_0\; ,\label{LSSH2}\\
u(t,x) &=& u(t,x+2\pi).\label{LSSH3}
\end{eqnarray}

Let $A=-(1+\partial_{xx})^{2}$. It is well-known that $A$
generates a compact analytic semigroup $e^{tA}$ in $C^{0}_{per}$.
Since
\begin{equation}\label{B_0}
B_{0}=\sum_{i=1}^{N^{'}}b^{2}_{i}<\infty,
\end{equation}
  the stochastic convolution $\int_{0}^{t}e^{A(t-\tau)}dW(x,t)$
is continuous for $(x,t)$; see \cite{PZ92}.

For convenience, we assume $N^{'}<\infty$. Let $\langle \cdot,
\cdot\rangle$ be the scalar product of $H=L^2_{per}(0, 2\pi)$,
$|\cdot|$ the norm, and $V=H^{1}\cap H$. Let
$$H_{l}=span\{e_{i},0<i \leq N\},\;\;\; H_{h}=span\{e_{i},i > N\},$$
which are two subspaces of $H$. We will call $H_l$ the low mode
space and $H_h$ the high mode space. Clearly, $H=H_{l}\bigoplus
H_{h}$. Let $P_{l}$ and $P_{h}$ be the orthogonal projections onto
$H_{l}$ and $H_{h}$, respectively. Denoting $u(t)=(l(t),h(t))$,
then (\ref{LSSH1})-(\ref{LSSH3}) is equivalent to the following
equations on $H$
\begin{eqnarray}
\dot{l}&=&P_{l}Al+\varrho l-P_{l}(l+h)^{3}+P_{l}F(x,t)\;,\;\;\;\;\;\;
l(0)=l_{0};\label{LM}\\
\dot{h}&=&P_{h}Ah+\varrho
h-P_{h}(l+h)^{3}\;,\;\;\;\;\;\;\;\;\;\;\;\;\;\;\;\;
 \;\;\;\;h(0)=h_{0},\label{HM}
\end{eqnarray}
where $l_{0}=P_{l}u(0)$ and $h_{0}=P_{h}u(0)$. For any given
$T>0$, since the noise is additive we can treat the equation
pathwise. Then the deterministic theory implies that
(\ref{LSSH1})-(\ref{LSSH3}) admits a mild solution with transition
probability $P(t,x,\Gamma)$ for $ \Gamma \in{\mathcal{B}}(H)$, the
Borel $\sigma$-algebra in $H$. Here, we treat
(\ref{LSSH1})-(\ref{LSSH3}) in another way that for $l_{0}\in
H_{l}$, $h_{0}\in H_{h}$, (\ref{LM})-(\ref{HM}) has a unique
solution $(l,h)\in X=C([0,T],H)$, where $l\in C(0,T;H_{l})$ and
$h\in C(0,T; H_{h})$. Notice that, for given $l(t)$ we can solve
(10) with $h_{0}\in H_{h}$. This solution is written as
$\Phi(l,h_{0})$. Much more such decomposition is discussed in
\cite{EL02}.

Let $l_{1}=l$, $ h_{1}=h$ for $t\in[0,T]$, $l_{k}$ and $h_{k}$ be
the solutions on $[0, kT]$, $l_{m,k}$ and $h_{m,k}$ be the
solutions on $[mT, kT]$ for $0\leq m\leq k$. Similar for $u_{m}$
and $u_{m,k}$.

For (\ref{LM}) and (\ref{HM}), we have the following proposition,
which displays the main idea of the paper.
\begin{prop}
For any $u_{0}^{1}$, $u_{0}^{2}\in H$, if
$l^{1}_{m,k}=l^{2}_{m,k}$ then there exists a constant $\lambda>0$
such that
$$|u^{1}-u^{2}|\leq e^{-\lambda (t-mT)}|u^{1}(mT)-u^{2}(mT)|
\;,\;\; t\in[mT, kT].$$
\end{prop}
\begin{proof}
Since $l^{1}_{m,k}=l^{2}_{m,k}$, we only consider the high mode.
Let $\rho=h^{1}-h^{2}$. Then $\rho$ satisfies
$$\dot{\rho}=A\rho+\varrho \rho-P_{h}((u^{1})^{3}-(u^{2})^{3}).$$
Taking the scalar product of the above equation with $\rho$ in
$H_{h}$ and notice that
$$\langle P_{h}((u^{1})^{3}-(u^{2})^{3}), \rho\rangle\geq 0,$$
we have
$$\partial_{t}|\rho|^{2}\leq (\alpha_{N}+\varrho)|\rho|^{2},$$
where $\alpha_{N}$ is the $Nth$ eigenvalue of $A$. If we choose
$N$ so large that $\alpha_{N}+\varrho<0$, then
$$|\rho|^2\leq e^{-\lambda (t-mT)}|\rho(mT)|\;,\;\; t\in[mT, kT],$$
for some constant $\lambda$. This completes the proof.
\end{proof}

\begin{remark}
This proposition implies that, if we want to get the distributions
of the LSSH system converged to the unique invariant measure, it
is enough to prove the following: for any two solutions $u^{1}$
and $u^{2}$, the probability of the event that those low modes
coincide is exponentially tends to $1$, as $t\rightarrow \infty$.
But, in fact, we do not need that all solutions satisfy the above
proposition. By the coupling result in Lemma 4.1 below, we just
need that the solutions that do not grow fast satisfy the
proposition. This will also be used for the NLSSH system.
\end{remark}

\section{Some Estimates}

In this section, we derive some estimates of the solutions for the
LSSH system, which will be used in the following. We will work on
a probability space $(\Om, \mathcal{F}, \mathbb{P})$ generated by
$\{W(t)\}$. We associate $\Om$ with the canonical space generated
by all $dw_i(t)$. Denote $\mathcal{F}$  as the associated
$\sigma$-algebra generated by $\{W(t)\}$ with $\mathbb{P}$   the
probability measure. Expectations with respect to $\mathbb{P}$
will be denoted by $\mathbf{E}$.\\

The following lemma describes the growth rate of $|u(t)|^2$.
\begin{lemma}
There is a positive constant $C_{1}$, such that for any $r>0$,
$${\mathbb{P}}\Big\{ |u(t)|^{2}+\int_0^t|u(s)|^2ds \leq|u(0)|^2+C_{1}t+r\;,
for\; t\geq 0 \Big\}\geq 1-e^{-r}.$$
\end{lemma}
\begin{proof}
Applying Ito's formula to $|u(t)|^{2}$, one yields
$$ \frac{1}{2}\frac{d}{dt}|u(t)|^{2}=\langle Au,u\rangle
+\langle \varrho u,u\rangle-\langle u^{3},u\rangle+\frac{1}{2}B_{0}
+\langle F(x,t),u\rangle.$$
Note that
\begin{eqnarray*}
\langle Au+\varrho u,u\rangle&=&\langle -(1-\varrho)u-2\partial_{xx}u
-\partial^{2}_{xx}u,u\rangle \\
&\leq&-(1-\varrho)|u|^{2}+|\partial_{xx}u|^2+|u|^{2}
-|\partial_{xx}u|^{2} \\
&=&\varrho |u|^{2}
 \end{eqnarray*}
and $$|u|\leq(2\pi)^{\frac{1}{4}}|u|_{L^{4}}.$$
Then we have
\begin{eqnarray*}
|u|^{2}+\int_0^t|u(s)|^2ds&\leq&
|u(0)|^{2}+\int_{0}^{t}[2(\varrho+\frac{1}{2})|u|^{2}-\frac{1}{2\pi}|u|^{4}+B_{0}]ds
+2\int_{0}^{t}\langle F,u\rangle\\
&=&
|u(0)|^{2}+\int_{0}^{t}[2(\varrho+b^{2}_{\max}+\frac{1}{2})|u|^{2}-\frac{1}{2\pi}
|u|^{4}+B_{0}]ds+M_{t}\\&&-\langle M\rangle_{t}/2
-(\int_{0}^{t}2b^{2}_{\max}|u|^{2}ds-\langle M\rangle_{t}/2)
\end{eqnarray*}
where $b_{\max}=\max_ib_i, M_{t}=2\int_{0}^{t}\langle F, u\rangle$
and $\langle M\rangle_{t}$ is the quadratic variation.

Note that
 $$\langle M\rangle_{t}/2=2\sum_{i}b^{2}_{i}\int_{0}^{t}u_{i}^{2}ds
 \leq 2b^{2}_{\max} \int_{0}^{t}|u|^{2}ds. $$
That is,
$$|u(t)|^{2}+\int_0^t|u(s)|^2ds\leq|u(0)|^{2}+C_{1}t+M_{t}-\langle M\rangle_{t}/2,$$
where
$C_{1}=\max_{x}\{2(\varrho+b^{2}_{\max}+\frac{1}{2})x^{2}-\frac{1}{2\pi}x^{4}+B_{0}\}.$

Then the classical supermartingale inequality implies
\begin{eqnarray*}
&&{\mathbb{P}}\Big\{|u(t)|^{2}+\int_0^t|u(s)|^2ds\leq |u(0)|^{2}+C_{1}t+r\Big\}\\
&\geq&  {\mathbb{P}}\Big\{M_{t}-\langle M\rangle_{t}/2\leq r\Big\}\\
&\geq& {\mathbb{P}}\Big\{\exp(M_{t}-\langle M\rangle_{t}/2)\leq e^{r}\Big\}\\
 &\geq& 1-e^{-r}.
\end{eqnarray*}
This completes the proof.
\end{proof}

We also need the following estimation of the mean value of the
solution.\\

\begin{lemma}
For any $t\geq 0$, we have
\begin{equation*}
{\mathbf{E}}|u(t)|^{2}\leq e^{-\alpha t}{\mathbf{E}}|u_{0}|^{2}+R,
\end{equation*}
where $\alpha$ and $R$ are some positive constants.
\end{lemma}
\begin{proof}
Applying It$\hat{o}$ formula to $|u|^{2}$ and
taking the mean value, we find
\begin{eqnarray*}
\textbf{E}|u|^{2}+\alpha\int_{0}^{t}\textbf{E}|u|^{2}ds&\leq&
\textbf{E}|u_{0}|^{2}+\int_{0}^{t}[(2\varrho-\alpha)\textbf{E}|u|^{2}-
\frac{1}{2\pi}\textbf{E}|u|^{4}+B_{0}]ds\\
&\leq& \textbf{E}|u_{0}|^{2}+R\alpha t,
\end{eqnarray*}
where $\alpha$ is some positive number and
 $R\alpha=\max_{x}\{(2\varrho-\alpha)x^{2}-\frac{1}{2\pi}x^{4}+B_{0}\}$.
Then Gronwall inequality yields the results.
\end{proof}

\section{Coupling Approach for Ergodicity}
We start this section with some notations and terminology. Let $H$
be a separable Hilbert space with $\sigma$-algebra
$\mathcal{B}(H)$ and let $\mathcal{M}(H)$ be the space of sigh
measures with bounded variation. We denote by $\mathcal{P}(H)$ the
set of probability measures $\mu\in\mathcal{M}(H)$. For
$\mu\in\mathcal{P}(H)$ we define the variation norm as
$$||\mu||_{var}=\sup_{A\in\mathcal{B}(H)}|\mu(A)|.$$
If $\mu_{1}$, $\mu_{2}\in\mathcal{P}(H)$ are absolutely continuous
with respect to a fixed Borel measure $m$, then we have
$$
||\mu_{1}-\mu_{2}||_{var}=\frac{1}{2}\int_{H}|\rho_{1}(u)-\rho_{2}(u)|dm(u),
$$
where $\rho_{i}$, $i=1,2,$ is the density of $\mu_{i}$ with respect
to $m$. For more see \cite{KS01}.\\

The coupling is a well-known effective tool for studying
finite-dimensional Markov chains \cite{Lin92}. To our knowledge,
\cite{Mu93} is the first paper  using a coupling approach for
invariant measures of the stochastic heat equation. The main idea
of coupling is to construct two random variables $\xi_{1}$ and
$\xi_{2}$, for   two measures $\mu_{1}$ and $\mu_{2}$ that we are
concerned with, and study the two measures through constructing
the two random variables. In this section, we will use the maximal
coupling, namely, we construct $\xi_{1}$ and $\xi_{2}$ such that
\begin{equation*}
||\mu_{1}-\mu_{2}||_{var}={\mathbb{P}}\Big\{\xi_{1}\neq\xi_{2}\Big\}.
\end{equation*}
About the maximal coupling, we refer to \cite{KS01} and \cite{Lin92}.\\

In this section, we always use the maximal coupling
$(l^{1}_{1},l_{1}^{2})$ for the distributions of $(P_{l}u^{1},
P_{l}u^{2})$ on $C(0,T;H_{l})$. By the discussion in Sec. 2,
we can define the coupling solution $u^{i}(t)$ as
\begin{equation*}
u^{i}(t)=\Big(l_{k}^{i}(t),\Phi(l^{i}_{k},
h^{i}((k-1)T))\Big),\;\; t\in[(k-1)T, kT],\;\;k\geq 1,\;\;i=1,2
\end{equation*}
for (5)-(7) on the product space $(\Omega^{k}, {\mathcal{F}}^{k},
{\mathbb{P}}^{k}) \times X^{k}$,
where
\begin{equation*}
\Omega^{k}=\Omega\times\cdots\times\Omega\;,\;{\mathcal{F}}^{k}
={\mathcal{F}}\times\cdots
\times{\mathcal{F}}\;,\;{\mathbb{P}}^{k}
={\mathbb{P}}\times\cdots\times{\mathbb{P}}.
\end{equation*}
and $X^{k}=X\times\cdots\times X$. From Theorem 4.2 in
\cite{KS01}, it follows that $u^{i}$ is measurable coupling
solutions of (\ref{LSSH1})-(\ref{LSSH3}). Notice that, the
coupling solutions are not the solutions of
(\ref{LSSH1})-(\ref{LSSH3}), but they have the same distributions
with the solutions of (\ref{LSSH1})-(\ref{LSSH3}), since we just
change $F(x,t)$ as another   noise with the same distribution.
Then, we will call the coupling solutions the weak solutions as in
\cite{KS03}. Obviously, the results in Sec. 2 still hold true if
$u(t)$ is a weak solution of (\ref{LSSH1})-(\ref{LSSH3}).\\

Now, for those weak solutions we define some function spaces which
play a role in our approach. This arises from Remark 2.2. These
spaces are also introduced in \cite{KS03}. For any given real
number $T>0$ and integer number $k\geq 1$, we define $S_{0}(m,k$)
as the set of functions $(u^{1}(t), u^{2}(t))$, $t\in [mT,kT]$,
$0\leq m\leq k$, such that
\begin{eqnarray}
u^{i}\in L^{2}(0,kT;V)\cap C(0,kT;H),\label{C1}\\
P_{l}u^{1}(t)=P_{l}u^{2}(t)\;,\;\;\;|u ^{1}(mT)|\vee|u^{2}(mT)|\leq D,
\label{C2}\\
{\mathcal{E}}^i(t,mT)\leq r+(C_{1}+1)(t-mT)\;, \;\;mT\leq t\leq
kT\;,\; i=1,2,\label{C3}
\end{eqnarray}
where $a\vee b=\max(a, b)$ for $a$, $b\in\R$ and
${\mathcal{E}}^i(t, mT)=|u^{i}(t)|^2+\int_{mT}^{t}|u^{i}(s)|^2ds$.
The set $S_{0}(m,k)$ is called the coupled set. Then the set
$$R(k)=(C(0,kT;H)\cap L^{2}(0,kT;V))\setminus \cup_{m=0}^{k}S_{0}(m,k)$$
should be the exponentially small set by Proposition 2.1. \\

Let
$$S(m,k)=S_{0}(m,k)\setminus S(m-1,k),\;\;0\leq m\leq k $$
where $S(-1,k)$ is the null set. In the following, we always
assume $|u^{1}_{0}|^{2}$ and $|u^{2}_{0}|^{2}$ have the finite
mean value. We observe that the orbit in $S(m,k-1)$ enters either
$S(m,k)$ or $R(k)\cup S(k,k)$ for $0\leq m\leq k-1$, and the orbit
in $S(m,k)$ comes from the set $R(m)\cup S(m,m)$; see Fig.1.
\medskip

{\small

 \setlength{\unitlength}{.4cm}
 \begin{picture}(40,22)
 \put(0,20){S(0,1)}
 \put(4,20){S(0,2)}
 \put(2.6,20.3){\vector(1,0){1.5}}
 \put(6.6,20.3){\vector(1,0){1.2}}
 \put(7.9,20){$\cdots$}
 \put(9.2,20.3){\vector(1,0){1.2}}

 \put(10.5,20){S(0,m)}
 \put(13.5,20.3){\vector(1,0){3}}
 \put(16.7,20){S(0,m+1)}
 \put(21,20.3){\vector(1,0){1.2}}
 \put(22.6,20){$\cdots$}
 \put(24.3,20.3){\vector(1,0){1.2}}

 \put(25.7,20){S(0,k-1)}
 \put(32.5,20){S(0,k)}
 \put(29.3,20.3){\vector(1,0){3}}
 \put(4,18){S(1,2)}
 \put(6.6,18.3){\vector(1,0){1.2}}
 \put(7.9,18){$\cdots$}
 \put(9.2,18.3){\vector(1,0){1.2}}

 \put(10.5,18){S(1,m)}
 \put(13.5,18.3){\vector(1,0){3}}
 \put(16.7,18){S(1,m+1)}
 \put(21,18.3){\vector(1,0){1.2}}
 \put(22.6,18){$\cdots$}
 \put(24.3,18.3){\vector(1,0){1.2}}

 \put(25.7,18){S(1,k-1)}
 \put(32.5,18){S(1,k)}
 \put(29.3,18.3){\vector(1,0){3}}

 \put(9.5,14){S(m-1,m)}
 \put(13.8,14.3){\vector(1,0){3}}
 \put(16.7,14){S(m-1,m+1)}

 \put(16.7,12){S(m,m+1)}
 \put(21.2,12.3){\vector(1,0){1.2}}
 \put(22.75,12){$\cdots$}
 \put(24.3,12.3){\vector(1,0){1.2}}

 \put(25.7,12){S(m,k-1)}
 \put(32.5,12){S(m,k)}
 \put(29.4,12.3){\vector(1,0){3}}

 \put(25,9){S(k-2,k-1)}
 \put(32.5,9){S(k-2,k)}
 \put(29.4,9.3){\vector(1,0){3}}

 \put(32.5,7){S(k-1,k)}

 \put(0,3.5){S(1,1)}
 \put(0,2.5){$\cup$}
 \put(1,2.5){R(1)}
 \put(1,4){\vector(1,4){3.5}}
 \put(3,3.3){\vector(1,0){1.2}}

 \put(4,3.5){S(2,2)}
 \put(4,2.5){$\cup$}
 \put(5,2.5){R(2)}
 \put(1.6,19.8){\vector(1,-4){3.9}}
 \put(7.5,3){$\cdots$}

 \put(9.5,3.5){S(m,m)}
 \put(9.5,2.5){$\cup$}
 \put(10.5,2.5){R(m)}
 \put(10.5,4){\vector(1,1){7.8}}
 \put(13.5,3){\vector(1,0){3}}

 \put(16.3,3.5){S(m+1,m+1)}
 \put(16.8,2.5){$\cup$}
 \put(17.8,2.5){R(m+1)}
 \put(23,3){$\cdots$}

 \put(25,3.5){S(k-1,k-1)}
 \put(25,2.5){$\cup$}
 \put(26,2.5){R(k-1)}
 \put(26.5,4){\vector(3,1){8}}
 \put(29,3){\vector(1,0){3}}

 \put(32.5,3.5){S(k,k)}
 \put(32.5,2.5){$\cup$}
 \put(33.5,2.5){R(k)}

 \put(29.3,20.3){\vector(1,-4){4}}
 \put(29.3,18.3){\line(1,-3){2}}
 \put(29.4,12.3){\line(1,-2){3.7}}
 \put(29.4,9.3){\line(1,-1){3}}

 \put(13.5,20.3){\vector(1,-4){4}}
 \put(13.5,18.3){\line(1,-3){2}}
 \put(13.8,14.3){\line(1,-2){2.4}}

 \put(33,16.5){$\cdot$}
 \put(33,15.5){$\cdot$}
 \put(33,14.5){$\cdot$}

 \put(26,16.5){$\cdot$}
 \put(26,15.5){$\cdot$}
 \put(26,14.5){$\cdot$}

 \put(19,16.5){$\cdot$}
 \put(19,15.5){$\cdot$}

 \put(15, 0){Figure 1: Orbit   transition}
 \end{picture}

\small}

\medskip

Then, what needs to be proved is that the event that $u^{1}$ and
$u^{2}$ are in $R(k)\cup S(k,k)$ has exponentially small
probability. Clearly, this demands that the events that the orbits
in $S(m,k-1)$ enter $R(k)\cup S(k,k)$ or do not enter $S(m,k)$
have exponentially small probability.\\

\begin{lemma}
Fix $T>0$ large enough. Let $\lambda_1$ and $\lambda_2$ be the
distributions of $l_{m,k}^{1}$ and $l_{m,k}^{2}$, respectively,
for $(u^{1}, u^{2})\in S(m,k-1)$, $0\leq m<k-1$. If (\ref{Mode})
holds and $N$ is large enough, then
$$||\lambda_{1}-\lambda_{2}||_{var}\leq e^{-\gamma (k-m)T}$$
for some positive constant $\gamma$.
\end{lemma}
\begin{proof}
Let $l^{1}$ and $l^{2}$ be the maximal coupling for $\lambda_{1}$ and
$\lambda_{2}$, respectively. There are two cases to be distinguished
for $t\in[0,\;T]$.\\

$\textbf{Case 1.} {\mathcal{E}}^i((k-1)T+t,mT)<r+(C_{1}+1)((k-1)T+t-mT),
\;i=1,2.$\\

Let
$$l_{0}^{i}=P_{l}u^{i}(mT),\;\;\;h_{0}^{i}=P_{h}u^{i}(mT),\;\;i=1,2.$$
Then $l^{1}_{0}=l^{2}_{0}=l$. Consider the LSSH system in low mode
space
\begin{equation*}
\frac{dl^{i}(t)}{dt}=[Al^{i}+\varrho
l^{i}-P_{l}(l^{i}+h^{i})^{3}]+P_{l}F(x,t)\;,\;\;i=1,2.
\end{equation*}
Let
$B(l,h_{0}^{1},h_{0}^{2})=P_{l}(l+h^{1})^{3}-P_{l}(l+h^{2})^{3}$.
Then
\begin{eqnarray*}
|B(l,h_{0}^{1},h_{0}^{2})|^{2}&\leq&\sup_{y\in L^{2}_{l}
 \; and\; |y|=1}|\langle
P_{l}(l+h^{1})^{3}-P_{l}(l+h^{2})^{3}, y\rangle|^{2}\\
&\leq&|(l+h^{1})^{3}-(l+h^{2})^{3}|^{2}\\
&\leq&|h^{1}-h^{2}|^{2}|(u^{1})^{2}+(u^{1})(u^{2})+(u^{2})^{2}|^{2}
\end{eqnarray*}
Proposition 2.1 tells us that
$$|h^{1}-h^{2}|^{2}\leq 4D^{2}e^{-\lambda((k-1)T+t-mT)}.$$
Then, there exists some positive constant $\lambda^{'}$, such that
$$|B(l,h_{0}^{1},h_{0}^{2})|^{2}\leq D^{2}e^{-\lambda^{'}(k-1-m)T},$$
since $|u^{i}(t)|^{2}$ at most increases polynomially. Thus
\begin{equation}
\int_{0}^{T}|B(l,h_{0}^{1},h_{0}^{2})|^{2} dt
\leq K_{0}e^{-\lambda^{'}(k-m-1)T}
\end{equation}
for some constant $K_{0}>0$. Now we apply the Girsanov's formula
to $\lambda_{i}$. Let
\begin{equation*}
\beta(t,\omega)=\textbf{b}^{-1}B(l,h_{0}^{1},h_{0}^{2}),
\end{equation*}
where $\textbf{b}$ is the $N\times N$ diagonal matrix with
diagonal elements $b_{j}$, $j=1, \cdots, N$.
Then, Girsanov's formula yields
\begin{equation}
\lambda_{1}(dl)=e^{G(l)}\lambda_{2}(dl),
\end{equation}
where
\begin{equation*}
G(l)=-\int_{0}^{T}(\beta,
\textbf{b}^{-1}dw_{t})-\frac{1}{2}\int_{0}^{T}|\beta|^{2}dt.
\end{equation*}

For (\ref{C2})
\begin{eqnarray*}
\textbf{E}\exp(2G(l))&=&\textbf{E}\exp\Big(-2\int_{0}^{T}(\beta,
\textbf{b}^{-1}dw_{t})-\int_{0}^{T}|\beta|^{2}dt\Big)\\
&\leq&\Big(\textbf{E}\exp(-4\int_{0}^{T}(\beta,
\textbf{b}^{-1}dw_{t})-8\int_{0}^{T}|\beta|^{2}dt)\Big)^{1/2}\\
&&\times\Big(\textbf{E}\exp(6\int_{0}^{T}|\beta|^{2}dt)\Big)^{1/2}\\
&&\leq e^{3K_{N}D^{2}e^{-\lambda^{'}(k-m-1)T}}
\end{eqnarray*}
with some constant $K_{N}$. Hence,
\begin{eqnarray*}
||\lambda_{1}-\lambda_{2}||_{var}&=&\frac{1}{2}\int_{C(0,T;H)}
\Big|1-\frac{d\lambda_{2}}{d\lambda_{1}}\Big|d\lambda_{1}\\
&\leq&
\frac{1}{2}\Big(\int_{C(0,T;H)}\Big|1-\frac{d\lambda_{2}}
{d\lambda_{1}}\Big|^{2}d\lambda_{1}\Big)^{1/2}\\
&\leq&\frac{1}{2}\Big(e^{3K_{N}D^{2}e^{-\lambda^{'}
(k-m-1)T}}-1\Big)^{1/2}\\
&\leq&  e^{-\gamma_{0}(k-m)T}
\end{eqnarray*}
for some positive constant $\gamma_{0}$, since $T>0$ is large enough. \\

$\textbf{Case 2.} {\mathcal{E}}^i((k-1)T+t,mT)\geq
r+(C_{1}+1)((k-1)T+t-mT)$ for some $t, i=1,2$.

This is the direct result of Lemma 3.1. In fact, the inequality
$${\mathcal{E}}^i((k-1)T+t,mT)\geq r+(C_{1}+1)((k-1)T+t-mT)$$
implies that
\begin{equation*}
{\mathcal{E}}^i((k-1)T+t,mT)\geq
|u(mT)|^{2}+C_{1}((k-1)T+t-mT)+(r-|u(mT)|^{2})+(k-m-1)T.
\end{equation*}
Taking $\gamma=\min\{1,\gamma_{0}\}$ completes the proof.
\end{proof}

\begin{remark}
Lemma 4.1 implies that, the probability of the event that the
orbits in $S(m,k-1)$, $0\leq m < k-1$ enter $R(k)\cup S(k,k)$
is small.
\end{remark}

To our aim we also need the following lemma.\\

\begin{lemma}
Let (\ref{Mode}) hold. If $(u^{1},u^{2})\in R(k-1)$ satisfies
 $|u^{1}((k-1)T)|\vee|u^{2}((k-1)T)|\leq D$, then
\begin{equation*}
{\mathbb{P}}\Big\{(u^{1},u^{2})\in S(k,k)\Big\}\geq d^{*}>0,
\end{equation*}
for some constant $d^{*}>0$ depending only on $D$.
\end{lemma}
\begin{proof}
From Lemma 3.1, we have
\begin{equation}\label{event}
{\mathbb{P}}\Big\{{\mathcal{E}}^i(t,(k-1)T)\leq
r+(C_{1}+1)(t-(k-1)T)\Big\}\geq c^{*}
\end{equation}
for some positive constant $c^{*}$ which depends only on $D$. For
the purpose of this lemma, we compare $l^{i}(t)$ with the
following standard diffusion process
\begin{equation}\label{Aidsys}
\dot{\tilde{l}}(t)=A\tilde{l}+(\varrho-\nu)\tilde{l}+P_{l}F(x,t).
\end{equation}
Here, $\nu$ is large enough such that (\ref{Aidsys}) is linearly
dissipative. Let $\tilde{l}^{i}$ be the solutions of
(\ref{Aidsys}) with the initial value $l^{i}((k-1)T)$. Let
$\lambda_{i}(t)$ and $\tilde{\lambda}_{i}(t)$ be the distribution
of $l^{i}_{k-1,k}$ and $\tilde{l}^{i}_{k-1,k}$ restricting on
$C((k-1)T, t; H_{l})$, respectively, for $i=1, 2$. Since
(\ref{Aidsys}) is linearly dissipative, we have
\begin{equation*}
\int \Big[\frac{d\tilde{\lambda}_{1}}{d\tilde{\lambda}_{2}}\Big]^{2}
 d\tilde{\lambda_{1}}\leq c_{1}^{*}
\end{equation*}
for some constant $c_{1}^{*}>0$ which depends only on the upper bound $D$.\\

By the similar method in the proof of  Lemma 4.1, we have
\begin{equation*}
\int
\Big[\frac{d\lambda_{i}}{d\tilde{\lambda_{i}}}\Big]^{2}d\tilde{\lambda}_{i}
\leq c_{2}^{*}\;\;i=1\;,\;2
\end{equation*}
for some positive constant $c_{2}^{*}$ which depends only on $D$.
Since at this time
$|B^{i}|=|\nu\tilde{l}^{i}+P_{l}(l^{i}+h^{i})^{3}|$ is bounded as
we consider the event of (\ref{event}) (See \cite{EL02}).  By
H$\ddot{o}$lder inequality, we can derive
\begin{equation*}
\int\Big[\frac{d\lambda_{1}}{d\lambda_{2}}\Big]^{2}d\lambda_{1}
 \leq c_{3}^{*}
\end{equation*}
for some positive constant $c_{3}^{*}$ which depends only on $D$.
Then by the Lemma C.1 of \cite{Mat02}, we have
\begin{equation*}
\int\Big|1\wedge\frac{d\lambda_{1}}{d\lambda_{2}}\Big|d\lambda_{2}
\geq c_{4}^{*}
\end{equation*}
for some constant $c_{4}^{*}>0$   depending only on $D$. Notice
that $||\lambda_{1}-\lambda_{2}||_{var}
=1-||\lambda_{1}\wedge\lambda_{2}||_{var}$, we have
\begin{equation*}
||\lambda_{1}-\lambda_{2}||_{var}\leq d^{*}=1-c_{4}^{*}.
\end{equation*}
Then ${\mathbb{P}}\Big\{(u^{1},u^{2})\in S(k,k)\Big\}\geq
1-d^{*}$. This completes the proof.
\end{proof}

Then, we have

\begin{prop}
Let (\ref{Mode}) hold. For $(u^{1},u^{2})\in R(k-1)\cup
S(k-1,k-1)$, there is a constant $0<\delta<1$, such that
\begin{equation*}
{\mathbb{P}}\Big\{(u^{1},u^{2})\in R(k) \Big\} \leq
\delta {\mathbb{P}}\Big\{ (u^{1},u^{2})\in R(k-1)\cup S(k-1,k-1)\Big\}
\end{equation*}
\end{prop}
\begin{proof}
This is the direct result of Lemma 4.3. We omit the proof.
\end{proof}

By Lemma 4.1 and Proposition 4.4, we have the following important
proposition, which implies the existence of a unique invariant measure.

\begin{prop}
Suppose (\ref{Mode}) holds. There is a constant $0<\kappa<1$ such
that for any initial value $u^{1}_{0}$, $u^{2}_{0}\in H$,
\begin{equation*}
{\mathbb{P}}\Big\{(u^{1},u^{2})\in R(k)\cup S(k,k)\Big\}\leq
K\kappa^{k}
\end{equation*}
where $K$ is a constant   depending only on the initial value.
\end{prop}
\begin{proof}
We will consider $R(k)$ and $S(k,k)$ respectively. First, we have
\begin{eqnarray*}
&&{\mathbb{P}}\Big\{(u^{1},u^{2})\in R(k)\Big\}\\
&\leq& e^{-\gamma kT}\sum_{m=0}^{k-2}e^{\gamma mT}
{\mathbb{P}}\Big\{(u^{1},u^{2})\in S(m,k-1)\Big\}+\\
&& \delta{\mathbb{P}}\Big\{(u^{1},u^{2})\in S(k-1,k-1)\cup R(k-1)\Big\}\\
&\leq&(2\delta e^{\gamma T}+1)e^{-\gamma kT}\sum_{m=0}^{k-3}e^{\gamma mT}
{\mathbb{P}}\Big\{(u^{1},u^{2})\in R(m)\Big\}+\\
&&(2\delta e^{\gamma T}+1)e^{-\gamma kT}\sum_{m=0}^{k-3}e^{\gamma mT}
{\mathbb{P}}\Big\{(u^{1},u^{2})\in S(m,m)\Big\}+\\
&&(e^{-2\gamma T}+\delta+\delta^{2}){\mathbb{P}}
\Big\{(u^{1},u^{2})\in R(k-2)\Big\}+\\
&&(2\delta^{2}+e^{-2\gamma T}){\mathbb{P}}
\Big\{(u^{1},u^{2})\in S(k-2,k-2)\Big\}\\
&\leq & f_{1}e^{-\gamma
kT}+f_{2}\delta^{k-1}+f_{3}\delta^{k}\\
&\leq& K \kappa ^{k},
\end{eqnarray*}
where $f_{i}$, $i=1, 2, 3$, increases on $k$ at most polynomially.

For $S(k,k)$, we have
\begin{eqnarray*}
&&{\mathbb{P}}\Big\{(u^{1},u^{2})\in S(k,k)\Big\}\\&\leq&
e^{-\gamma kT}\sum_{m=0}^{k-2}e^{\gamma mT}
{\mathbb{P}}\Big\{(u^{1},u^{2})\in S(m,k-1)\Big\}+\\
&&\delta{\mathbb{P}}\Big\{(u^{1},u^{2})\in
S(k-1,k-1)\Big\}+{\mathbb{P}}\Big\{(u^{1},u^{2})\in
R(k-1)\Big\}.
\end{eqnarray*}
By the same analysis for $R(k)$, we also have
\begin{equation*}
{\mathbb{P}}\Big\{(u^{1},u^{2})\in S(k,k)\Big\}\leq K \kappa^{k}.
\end{equation*}
This completes the proof.
\end{proof}


\section{Ergodicity of the Local System}

In this section, we prove the existence and uniqueness of invariant
measure, and then the ergodicity in the LSSH system.\\

Let ${\mathcal{P}}(H)$ be the metric space of all probability
measures on $(H,{\mathcal{B}}(H))$ endowed with the metric
$||\cdot||^{*}_{L}$ defined by
\begin{equation*}
||\mu-\nu||^{*}_{L}=\sup\Big\{\int_{H}f d(\mu-\nu): ||f||_{L}\leq 1
\Big\},\;\;\;\mu,\;\nu\in{\mathcal{P}}(H)
\end{equation*}
where $f$ is a measurable function on $H$ and
\begin{equation*}
||f||_{L}=\sup_{x}|f(x)|+\sup\{\frac{f(x)-f(y)}{|x-y|}: x\neq y\in H\}.
\end{equation*}\\

It is well known that $||\cdot||^{*}_{L}$ generates the weak
topology and ${\mathcal{P}}(H)$ is a complete space under
$||\cdot||^{*}_{L}$. We define the semigroup $S(t)$ on
${\mathcal{P}}(H)$ generated by (\ref{LSSH1})-(\ref{LSSH3}) as
\begin{equation*}
S(t)\mu(\Gamma)=\int_{H}P(t,x,\Gamma)\mu(dx)
\;,\;\;\Gamma\in{\mathcal{B}}(H).
\end{equation*}
A measure $\mu$ on $H$ is an invariant measure if
\begin{equation*}
\mu(\Gamma)=S(t)\mu(\Gamma)=\int_{H}P(t,x,\Gamma)\mu(dx)
\;,\;\;\Gamma\in{\mathcal{B}}(H).
\end{equation*}
\\

We restrict $S(t)$ on the set ${\mathcal{P}}_{2}(H)$,
which is defined as
\begin{equation*}
{\mathcal{P}}_{2}(H)=\Big\{\mu\in{\mathcal{P}}(H):
\int_{H}|z|^{2}\mu(dz)<\infty\Big\}.
\end{equation*}
Then we have the following proposition.\\

\begin{prop}
S(t) maps ${\mathcal{P}}_{2}(H)$ into ${\mathcal{P}}_{2}(H)$.
\end{prop}
\begin{proof}
For any $\mu\in{\mathcal{P}}_{2}(H)$ and $t>0$, Lemma 3.2 yields
\begin{eqnarray*}
\int_{H}|y|^{2}S(t)\mu(dy)
&=&\int_{H}\Big(\int_{H}|y|^{2}P(t,z,dy)\Big)\mu(dz)\\
&=&\int_{H}{\mathbf{E}}|u(t,z)|^{2}\mu(dz)\\
&\leq&e^{-\alpha t}\int_{H}|z|^{2}\mu(dz)+R \\
&<&\infty.
\end{eqnarray*}
This completes the proof.
\end{proof}

Now we can draw the following main results.\\

\begin{theorem} \label{contraction}
If (\ref{Mode}) and (\ref{B_0}) hold with $N$ large enough, there
is a positive constant $\chi\in(0,1)$ such that for any $\mu_{1}$,
$\mu_{2}\in{\mathcal{P}}_{2}(H)$,
\begin{equation*}
||S(t)\mu_{1}-S(t)\mu_{2}||^{*}_{L}\leq
M(\mu_{1},\mu_{2})\chi^{t}\;, t\geq 0.
\end{equation*}
Here $M(\mu_{1},\mu_{2})$ only depends on $\mu_{1}$ and $\mu_{2}$.
\end{theorem}
\begin{proof}
We follow the proof of \cite{KS03}. Arbitrarily fixed $t>0$ and
let $k=k(t)$ be the smallest integer such that $t\leq k(t)T$,
where $T$ is the constant in Sec. 4, and let $u^{i}_{0}$, $i=1,2$,
be random variables in $H$ with distribution $\mu_{i}$. Let
$u^{1}(t)$ and $u^{2}(t)$ be the weak solutions of
(\ref{LSSH1})-(\ref{LSSH3}) on [0,kT]. Then we just need to show
that
\begin{equation}
p(t)={\mathbb{P}}\Big\{|u^{1}(t)-u^{2}(t)|>C_{1}e^{-\sigma_{1}
t}\Big\}\leq C_{2}e^{-\sigma_{2}t}
\end{equation}
where $C_{i}$, $i=1,2$, are positive constants only depending on
the initial functions and $\sigma_{i}$, $i=1,2$, are some positive
constants.\\

We consider the coupling solution defined in Sec. 4. Define the
following event
\begin{equation*}
G(k)=\Big\{(u_{k}^{1},u_{k}^{2})\in\bigcup_{m=0}^{[ck]}S(m,k)\Big\}
\end{equation*}
where $0<c<1$ is some constant.
Clearly,
\begin{equation}
p(t)\leq
{\mathbb{P}}\Big(G(k)^{c}\Big)+{\mathbb{P}}\Big(G(k)\cap\{|u^{1}-u^{2}|
>C_{1}e^{-\sigma_{1}t}\}\Big).
\end{equation}
We claim that
\begin{eqnarray}
&&{\mathbb{P}}\Big(G(k)^{c}\Big)\leq
C_{2}e^{-\sigma_{2}t}\;,\\
&&{\mathbb{P}}\Big(G(k)\cap \{|u^{1}-u^{2}|>C_{1}e^{-\sigma_{1}t}\}
\Big)=0.
\end{eqnarray}

First, we prove (24). Since all the orbits in $S(m,k)$ come from
$R(m)\cup S(m,m)$, we have
\begin{equation*}
G(k)^{c}\subset\bigcup_{m=[ck]+1}^{k}S(m,k)\cup R(k)\subset
\bigcup_{m=[ck]+1}^{k} R(m)\cup S(m,m).
\end{equation*}
Now by Proposition 4.5, we get
\begin{eqnarray*}
{\mathbb{P}}\Big(G(k)^{c}\Big)
&\leq&\sum_{m=[ck]+1}^{k}K\kappa^{m}
\leq K\kappa^{ck}(k-[ck])\\
&\leq& \tilde{K}\tilde{\kappa}^{k}
\end{eqnarray*}
for some proper $\tilde{\kappa}\in(0,1)$ and $\tilde{K}$.
Since $k\geq \frac{t}{T}$, we have
\begin{equation*}
{\mathbb{P}}\Big(G(k)^{c}\Big)\leq C_{2}e^{-\sigma_{2}t}
\end{equation*}
where $\sigma_{2}=-T^{-1}\ln\tilde{\kappa} $ and
$C_{2}=\tilde{K}$.\\

Next we prove (25). It is enough to prove that if
$(u^{1},u^{2})\in G(k)$ and $C_{1}$ is large enough
\begin{equation}
|u^{1}-u^{2}|\leq C_{1}e^{-\sigma t}
\end{equation}
Indeed, by the definition of $G(k)$ for $(u^{1},u^{2})\in G(k)$,
there is an integer $l$, $0\leq l\leq[ck]$, such that
$(u^{1},u^{2})\in S(l,k)$. Therefore, the relations
(\ref{C1})-(\ref{C3}) are satisfied. Then as Proposition 2.1, for
$u=u_{1}-u_{2},$ we have
\begin{equation*}
|u(t)|=|\rho(t)|\leq 2de^{-\lambda(t-lT)}.
\end{equation*}
Notice that $lT\leq ckT\leq c(t+T)$, then $t-lT\geq (1-c)t-cT$.
Hence, $|u(t)|\leq 2de^{\lambda cT}e^{-\lambda(1-c)t}$. Thus (25)
holds with $C_{1}=2de^{c\lambda T}$ and $\sigma_1=\lambda(1-c)$.

Finally, taking $M(\mu_{1},\mu_{2})=C_2$ and
$\chi=e^{-\sigma_{2}}$, we obtain the conclusion and the proof is
completed.
\end{proof}

\begin{coro}\label{ergodicity}
If enough modes are forced and the condition (\ref{Mode}) holds,
then the LSSH system has a unique invariant measure $\mu_{0}$,
such that for any $u\in H$ and $t>0$,
\begin{equation*}
||P(t,u,\cdot)-\mu_{0}||_{L}^{*}\leq M(|u|^{2})\chi^{t},
\end{equation*}
where $M(|u|^{2})$ is a positive constant and $\chi$ is as in Theorem 5.2.
\end{coro}


\section{Ergodicity of the Nonlocal System}

In this section, we turn to the following NLSSH system on the real line
\begin{eqnarray}
u_t-\varrho u+(1+\partial_{xx})^2u+uG*u^2&=&\tilde{F}(x,t)\;,\label{NLSSH1}\\
u(0) &=& u_0\; ,\label{NLSSH2}\\
u(t,x) &=& u(t,x+2\pi),\label{NLSSH3}
\end{eqnarray}
where
\begin{equation*}
 G*u^2=\int_0^{2\pi}G(x-\xi)u^2(\xi,t)d\xi.
\end{equation*}

We always assume that
\begin{equation}\label{Mode2}
\tilde{b}_{i}\neq 0\;,\;\;1\leq i\leq \tilde{N}\leq N''
\end{equation}
\begin{equation}\label{B0}
 \tilde{B}_0=\sum_{i=1}^{N''}\tilde{b}^2_i<\infty.
\end{equation}

First, we consider the following case: for every $x\in \R$
\begin{equation}\label{Ker1}
0<b\leq G(x)\leq a,
\end{equation}
where $a$ and $b$ are the positive constants (i.e., $G>0$ is a
positive kernel). Then we have the same energy estimation as LSSH
system.

\begin{lemma}
There is a positive constant $\tilde{C}_{1}$, such that for any
$r>0$,
$${\mathbb{P}}\Big\{ |u(t)|^{2}+\int_0^t|u(s)|^2ds \leq|u(0)|^2
+\tilde{C}_{1}t+r\;,for\; t\geq 0 \Big\}\geq 1-e^{-r}.$$
\end{lemma}
\begin{proof}
We only estimate the nonlinear term. In fact,
\begin{equation*}
\langle uG*u^2, u\rangle=\langle G*u^2, u^2\rangle\geq b|u|^4.
\end{equation*}
Set
$\tilde{C}_1=\max_x\Big\{{2(\varrho+\tilde{b}^2_{max}
+\frac{1}{2})x^2-bx^4+\tilde{B}_0}\Big\},$
where
\begin{equation}
   \tilde{b}_{max}=\max_{i}\tilde{b}_i. \label{bmax}
\end{equation}
By using the
same analysis as in Lemma 3.1, we can finish the proof.
\end{proof}

Note that, for the nonlocal nonlinearity, only the part of results
in Proposition 2.1 still holds. But this does not affect the
coupling result. Taking the same decomposition of NLSSH system as
one of LSSH system, we have

\begin{prop}
For any $u_{0}^{1}$, $u_{0}^{2}\in H$, if
$l^{1}_{m,k}=l^{2}_{m,k}$ and $u^1$, $u^2$ satisfy (\ref{C3}),
then there exists a constant $\tilde{\lambda}>0$ such that
$$|u^{1}-u^{2}|\leq e^{-\tilde{\lambda}(t-mT)}|u^{1}(mT)-u^{2}(mT)|
\;,\;\; t\in[mT, kT].$$
\end{prop}
\begin{proof}
Since $l^{1}_{m,k}=l^{2}_{m,k}$, we only consider the high mode.
Let $\rho=h^{1}-h^{2}$. Then, $\rho$ satisfies
$$\dot{\rho}=A\rho+\varrho \rho-P_{h}(u^1G*(u^{1})^{2}-u^2g*(u^{2})^{2}).$$
Taking the scalar product of the above equation with $\rho$ in
$H_{h}$ and notice that
\begin{eqnarray*}
&&\langle u^1G*(u^{1})^{2}-u^2G*(u^{2})^{2}, \rho\rangle \\
&=&\langle u^1G*(u^1)^2-u^1G*(u^2)^2+u^1G*(u^2)^2-u^2G*(u^2)^2,
\rho\rangle\\
&=&\langle u^1G*((u^1)^2-(u^2)^2),
\rho\rangle+\langle(u^1-u^2)G*(u^2)^2, \rho\rangle \\
&\geq&\langle u^1G*(u^1-u^2)(u^1+u^2), \rho\rangle,
\end{eqnarray*}
we have
\begin{eqnarray*}
\partial_{t}|\rho|^{2}&\leq& (\alpha _{\tilde{N}}+\varrho)|\rho|^{2}+|\langle
u^1G*(u^1-u^2)(u^1+u^2), \rho\rangle|,\\
&\leq&(\alpha
_{\tilde{N}}+\varrho)|\rho|^{2}+a(\frac{3}{2}|u^1|^2+|u^2|^2)|\rho|^2,
 \end{eqnarray*}
where $\alpha_{\tilde{N}}$ is the $\tilde{N}th$ eigenvalue of $A$.
Then
\begin{eqnarray*}
|\rho|^2&\leq&|\rho(mT)|^2\exp\Big\{\alpha_{\tilde{N}}+\varrho+
   a\int_{mT}^t(\frac{3}{2}|u^1|^2+|u^2|^2)ds\Big\}\\
   &\leq&|\rho(mT)|^2\exp\Big\{\Big(\alpha_{\tilde{N}}+\varrho+
   \frac{a}{t-mT}\int_{mT}^t(\frac{3}{2}|u^1|^2
   +|u^2|^2)ds\Big)\Big(t-mT\Big)\Big\}\\
   &\leq&|\rho(mT)|^2\exp\Big\{(\alpha_{\tilde{N}}+\varrho
   +\frac{5}{2}a\tilde{C}_1+a)(t-mT)\Big\},
\end{eqnarray*}
since we can always take $r<T$. If we choose $\tilde{N}$ so large
that $\alpha_{\tilde{N}}+\varrho+\frac{5}{2}a\tilde{C}_1+a<0$, then
$$|\rho|^2\leq e^{-\tilde{\lambda}(t-mT)}|\rho(mT)|\;,\;\; t\in[mT, kT],$$
for some constant $\tilde{\lambda}$. This completes the proof.
\end{proof}
\begin{remark}
The above proposition is weaker than Proposition 2.1, but it is
enough for the coupling result in Lemma 4.1.
\end{remark}

We have to give more estimation for the coupling result. For Lemma
4.1, we have to estimate
$|\tilde{B}(l,h_0^1,h_0^2)|=|P_lu^1G*(u^1)^2-P_lu^2G*(u^2)^2|$.
In fact,
\begin{eqnarray*}
&&|u^1G*(u^1)^2-u^2G*(u^2)^2|\\
&=&|(l+h^1)G*(l+h^1)^2-(l+h^2)G*(l+h^2)^2|\\
&\leq&|lG*[(l+h^1)^2-(l+h^2)^2]|\\
&&+|h^1G*(l^2+(h^1)^2+2lh^1)-h^2G*(l^2+(h^2)^2+2lh^2)|\\
&\leq&|2lG*[l(h^1-h^2)]|+|lG*[(h^1-h^2)(h^1+h^2)]|\\
&&+|(h^1-h^2)G*l^2|+|h^1G*(h^1)^2-h^2G*(h^2)^2|\\
&&+2|h^1G*(lh^1)-h^2G*(lh^2)|.
\end{eqnarray*}

From (\ref{Ker1}) we have
\begin{equation}
|2lG*[l(h^1-h^2)]|\leq 2a|h^1-h^2||l|^2.
\end{equation}
\begin{equation}
|lG*[(h^1-h^2)(h^1+h^2)]|\leq a|h^1-h^2||h^1+h^2||l|.
\end{equation}
\begin{equation}
|(h^1-h^2)G*l^2|\leq a|h^1-h^2||l^2|.
\end{equation}
\begin{eqnarray}
&&|h^1G*(h^1)^2-h^2G*(h^2)^2|\\
&\leq&|h^1G*(h^1)^2-h^2G*(h^1)^2|+|h^2G*(h^1)^2-h^2G*(h^2)^2|\\
&\leq& a|h^1-h^2||h^1|^2+a|h^1-h^2||h^1+h^2||h^2|.
\end{eqnarray}
Similarly, we have
\begin{equation}
 2|h^1G*(lh^1)-h^2G*(lh^2)|\leq 2a|h^1-h^2||h^1+h^2||l|.
\end{equation}

Then, from (33)-(39) and $|u|^2$ increases polynomially, we can
derive the coupling result in Lemma 4.1 by Proposition 6.2.

For Lemma 4.3, we have to estimate
$|\tilde{B}|=|\nu\tilde{l}+P_luG*(l+h)^2|$. Since
\begin{equation}
|uG*u^2|\leq a|u|^3,
\end{equation}
$|\tilde{B}|$ is bounded.\\

With all the above analysis, we draw the following conclusion.
\begin{theorem}  \label{im}
If (\ref{Mode2}) and (\ref{B0}) hold with $\tilde{N}$ large
enough, the NLSSH system with the positive kernel (\ref{Ker1}) has
a unique invariant measure, and it is exponentially attracted in
${\mathcal{P}}_{2}(H)$.
\end{theorem}

The above analysis for the positive kernel (\ref{Ker1}) does not
work for the nonnegative kernels. But we can still work for a
special nonnegative kernel. See \cite{LGDE00}.

Define
\begin{eqnarray}
 J(x)=\left\{
  \begin{array}{c l}
    c\exp(-\frac{1}{1-x^2}) & \mbox{if $x<1$},\\
    0 & \mbox{if  $x\geq 1$}.
  \end{array}
\right.
\end{eqnarray}
where
\begin{equation*}
c=\Big(\int_0^1\exp(-\frac{1}{1-x^2})dx\Big)^{-1}.
\end{equation*}
We further define that for $\delta>0$,
\begin{equation*}
J_{\delta}(x)=\delta^{-2}J(\frac{x}{\delta}).
\end{equation*}
Let $C_{0}(\Bar{I})$ be the space of continuous functions with
compact support in $I$. It is known that \cite{Fri}
\begin{equation}
||J_{\delta}*f-f||_{C_{0}(\bar{I})}\rightarrow 0 \;\; as \;\;
\delta\rightarrow 0,
\end{equation}
for any $f\in C_{0}(\Bar{I})$. Thus for any $\epsilon>0$, there is
a $\delta_{0}=\delta_{0}(\epsilon)>0$, such that
\begin{equation}\label{Ineq}
J_{\delta_0}*f\geq f-\epsilon.
\end{equation}\\

We consider a special kernel $G(x)=J_{\delta_{0}}(x)$, with
$\delta_{0}=\delta_{0}(\epsilon)$ as in (\ref{Ineq}), for the
NLSSH (\ref{NLSSH}). Then, $G$ satisfies
\begin{equation}
0\leq G\leq  \frac{c}{\delta_{0}^2}.
\end{equation}
Note that the lower bound of $G$ is $b=0$, so some estimation above
does not apply. But due to (\ref{Ineq}) we can still get the
similar energy estimates. And notice that we only use the lower
bound in the energy estimates, so we can have the same ergodic
result.

From \cite{Henry}, \cite{Pazy} or \cite{LGDE00}, we know that the problem
(\ref{NLSSH1})-(\ref{NLSSH3}) has a unique solution
$u(t,x,u_{0})\in C([0,T]; H)\cap L^{\infty}((0,T), H_0^2(I))$ and
\begin{equation*}
u(t,x,u_0)\in C_0(\Bar{I}),  \forall t\geq t_0>0.
\end{equation*}
Thus (\ref{Ineq}) implies that
\begin{equation}
G(x)*u^2\geq u^2-\epsilon.
\end{equation}
Then the estimates of Lemma 6.1 still holds with
\begin{equation}
\tilde{C}_1=\max_x\Big\{{2(\mu+\tilde{b}^2_{max}+\frac{1}{2}+\epsilon)x^2
-\frac{1}{2\pi}x^4+\tilde{B}_0}\Big\}.
\end{equation}

The other estimates does not depends on the lower bound of $G$.
Then we can derive the following result.

\begin{theorem}  \label{im2}
If (\ref{Mode2}) and (\ref{B0}) hold with $\tilde{N}$ large
enough, the NLSSH system with the nonnegative kernel
$J_{\delta_0}(x)$ has a unique invariant measure, and it is
exponentially attracted in ${\mathcal{P}}_{2}(H)$.
\end{theorem}

\section{Dynamical Difference between Local and Nonlocal Systems}

With the results on ergodicity of the local and nonlocal
stochastic Swift-Hohenberg equations in previous sections,
we can now compare some dynamical behavior.\\

From Proposition 2.1, we see that  the local stochastic
Swift-Hohenberg system is ergodic provided $\alpha_N+\varrho<0$.
Note that $\alpha_N=-(N^2-1)^2$. So we only need the number $N$ of
randomly forced modes to satisfy
\begin{equation*}
N>\sqrt{\sqrt{\varrho}+1}.
\end{equation*}
Note that the number $N$ does not   depend on the random forcing
term, i.e., it does not depend on the coefficients $b_i's$ in the
random forcing in (\ref{Noise1}).\\

However, from Proposition 6.2 the nonlocal stochastic
Swift-Hohenberg system with positive kernel ($0<b \leq G(x) \leq
a$) is ergodic provided that $\tilde{N}$ satisfies
$\alpha_{\tilde{N}}+\varrho+\frac{5}{2}a\tilde{C}_1+a<0$, where
$\tilde{C}_1=\tilde{B}_0+\frac{(\varrho+\tilde{b}^2_{max}+\frac{1}{2})^2}{b}$,
with $\tilde{B}_0$ defined in (\ref{B0}) and $\tilde{b}_{max}$
defined in (\ref{bmax}). Thus we only need the number of randomly
forced modes to satisfy
\begin{equation*}
\tilde{N}>\sqrt{\sqrt{\varrho+a+\frac{5}{2}a(\tilde{B}_0+
 \frac{(\varrho+\tilde{b}^2_{max}+\frac{1}{2})^2}{b})}+1}.
\end{equation*}
For the nonnegative kernel $G(x)$ with upper bound $a$, we only
need
\begin{equation*}
\tilde{N}>\sqrt{\sqrt{\varrho+a+\frac{5}{2}a(\tilde{B}_0+
2\pi(\varrho+\tilde{b}^2_{max}+\frac{1}{2}+\epsilon)^2)}+1}.
\end{equation*}
\\

Clearly, we have the following comparison between the local and
nonlocal stochastic Swift-Hohenberg models: (i) The number of
Fourier modes to be randomly excited for ensuring ergodicity of the
{\em local} stochastic Swift-Hohenberg system depends only on the
Rayleigh number through $\varrho>0$, which measures the difference
of the Rayleigh number from its critical convection onset value;
(ii) The number of Fourier modes to be randomly excited for
ensuring ergodicity of the {\em nonlocal} stochastic
Swift-Hohenberg system depends on the   bound of the kernel $G$ in
the nonlocal term, and the random term itself, as well as on
the Rayleigh number.

\bigskip

{\bf Acknowledgement.} We thank Dirk Blomker and Martin Hairer for
useful comments and discussions.

\end{document}